\documentclass[11pt,titlepage]{article}

\usepackage[margin=2cm]{geometry}
\usepackage{fancyhdr}
\usepackage{amssymb,amsmath}
\usepackage{float}

\usepackage{eucal,dsfont,accents,bbm,amsfonts,amsthm,amsopn,mathrsfs, multirow}
\usepackage{enumerate}
\usepackage[usenames]{color}
\usepackage{makeidx}

\usepackage{array}
\newcolumntype{C}{>{$}m{0.65cm}<{$}}

\newcommand\lm{\lambda}

\newtheorem{main}{Theorem}
\newtheorem{thm}{Theorem}

\newtheorem{defi}[thm]{Definition}
\newtheorem{lem}[thm]{Lemma}
\newtheorem{cor}[thm]{Corollary}

\newcommand\pf{\noindent{\em Proof.}~~}

\newcommand\lspan[1]{{\langle{#1}\rangle}}

\newcommand\agen[1]{{\langle\!\langle{#1}\rangle\!\rangle}}

\newcommand\ggen[1]{{{<}{#1}{>}}}

\title{Majorana Algebra for the Hoffman--Singleton Graph}
\author{Andries E. Brouwer and Alexander A. Ivanov\footnote{Imperial College London, London, UK, a.ivanov@ic.ac.uk }\,\footnote{Three Gorges Mathematical Research Center, Yichang, Hubei, China}\,\footnote{Institute for System Analysis ERC, CSC RAS, Moscow, Russia}}
\date{\today}

\begin{document}
	\maketitle
\begin{abstract}
Majorana theory is an axiomatic tool for studying the Monster group $M$ and its subgroups through the 196\,884-dimensional Conway--Griess--Norton algebra. The theory was introduced by A.~A.~Ivanov in 2009 and since then it experienced a remarkable development including the classification of Majorana representations for small (and not so small) groups. The group $U_3(5)$ is (isomorphic to) the socle of the centralizer in $M$ of a subgroup of order 25. The involutions of this $U_3(5)$-subgroup are $2A$-involutions in the Monster. Therefore, $U_3(5)$ possesses a Majorana representation (based on the embedding in the Monster). We prove that this is the unique Majorana representation of $U_3(5)$, calculate its dimension, which is $798$, and obtain a description in terms of the Hoffman--Singleton graph of which the automorphism group has $U_3(5)$ as an index 2 subgroup.
\end{abstract}	

\section{The Monster and Majorana}
The Monster group, which is the largest and most famous sporadic simple group,
has its minimal complex representation of dimension 196\,883. This number
was noticed by J.~McKay in the 1970's to be one less than the linear
coefficient of the modular invariant $J(q)$,
% standard notation: j(z) = q^{-1}+744+196884q+... for q = exp(2 \pi i z)
making the start of the Monstrous Moonshine. The underlying vector space, as shown by S.~P.~Norton in the early 1970's, carries invariant inner and algebra products. This algebra was more explicitly described by R.~L.~Griess \cite{g82} when constructing the Monster in 1980. Later the construction was reviewed and improved in various ways by J.~H.~Conway \cite{c84}. In particular Conway adjoined to the algebra an identity and defined $2A$-axes associated with $2A$-involutions in the Monster. The $2A$-axis $a(t)$ associated with an involution $t$ is an idempotent in the 2-dimensional space
$$C_{V_{M}}(C_M(t)),{\rm~where~}C_M(t) \cong 2 \cdot BM,$$
$BM$ standing for the Baby Monster sporadic simple group. Among other idempotents in this 2-space the $2A$-axis is characterized by having 1-eigenvector with multiplicity 1. The 196\,884-dimensional algebra is called the Conway--Griess--Norton algebra or simply the Monster algebra and is denoted by $V_M$.

\medskip\noindent
In \cite{iv09} essential properties of $2A$-axes in the Monster algebra were axiomatised in the following way.

\medskip\noindent
	Let $V$ be a real vector space equipped with a bilinear form $(~,~)$ and an algebra product $~\cdot~$. 
	
	\begin{itemize}
		\item[{\rm (M1)}]  $(~,~)$ is a symmetric positive-definite bilinear form on 
		$V$ that associates with $\cdot$ in the sense that 
		$(u,v\cdot w)=(u\cdot v,w)$ for all $u,v,w \in V$, and $\cdot$
		is a bilinear commutative non-associative
		algebra product on $V$;
		\item[{\rm (M2)}] the Norton inequality holds, so that 
		$(u \cdot u,v \cdot v)\ge (u\cdot v,u\cdot v)$
		for all $u,v \in V$.
	\end{itemize}
	
	\medskip\noindent
	A vector $a \in V$ is said to be a {\em Majorana axis} if it satisfies the following five conditions (M3) to (M7), where ${\rm ad}(a): v \mapsto a \cdot v$ is the adjoint operator of $a$ on $V$.
	\begin{itemize}	
		\item[{\rm (M3)}] $(a,a)=1$ and $a\cdot a=a$, so that $a$ is an idempotent of length 1;
		\item[{\rm (M4)}] ${\rm ad}(a)$ is semi-simple with spectrum $Sp=\{1,0,\frac{1}{4},\frac{1}{32}\}$:
		$$V=V_1^{(a)} \oplus V_0^{(a)} \oplus V_{\frac{1}{4}}^{(a)}
		\oplus V_{\frac{1}{32}}^{(a)},$$
		where $V_\mu^{(a)}=\{v \mid v \in V, a \cdot v=\mu v\}$ is the set of $\mu$-eigenvectors
		of ${\rm ad}(a)$ on $V$;
		\item[{\rm (M5)}] $V_1^{(a)}=\{\lm a \mid \lm \in {\mathbb R}\}$;
		
		\item[{\rm (M6)}] the linear transformation $\tau(a)$ of $V$ defined via
		$$\tau(a): u \mapsto (-1)^{32\mu}u$$
		for $u \in V^{(a)}_\mu$ with $\mu=1,0,\frac{1}{4},\frac{1}{32}$, preserves
		the algebra product ({\em i.e.} $u^{\tau(a)} \cdot v^{\tau(a)}=(u \cdot v)^{\tau(a)}$
		for all $u,v \in V$). The automorphism $\tau(a)$ is called the {\em Majorana involution} associated with the Majorana axis $a$;
		\item[{\rm (M7)}] if $V_+^{(a)}$ is the centralizer of $\tau(a)$ in $V$,
		so that $V_+^{(a)}=V_1^{(a)} \oplus V_0^{(a)} \oplus V_{\frac{1}{4}}^{(a)}$,
		then the linear transformation $\sigma(a)$ of $V_+^{(a)}$ defined via 
		$$\sigma(a): u \mapsto (-1)^{4\mu}u$$
		for $u \in V^{(a)}_\mu$ with $\mu=1,0,\frac{1}{4}$ preserves the restriction
		of the algebra product to the subalgebra $V^{(a)}_+$.
	\end{itemize}

	\medskip
	\noindent
	The conditions (M1) to (M7) imply that the eigenspaces $V_\mu^{(a)}$ of the
	adjoint action of $a$ satisfy the fusion rules described in the following table:
	
	\begin{center}
		\vspace{0.5cm}		
		\begin{tabular}{|c|cccc|}
			\hline
			&&&&\\
			$Sp$ & ~~~$1$~~~  &  $0$  & $\frac{1}{4}$ & $\frac{1}{32}$  \\
			&&&&\\
			\hline
			&&&&\\
			~~~$1$~~~ & $1$ & ~~$0$~~ & ~~~~$\frac{1}{4}$~~~~ & $\frac{1}{32}$  \\
			&&&&\\
			$0$ & $0$ & $0$ & $\frac{1}{4}$ & $\frac{1}{32}$  \\
			&&&&\\
			$\frac{1}{4}$  &  $\frac{1}{4}$  &  $\frac{1}{4}$ & $1,0$ & $\frac{1}{32}$\\
			&&&&\\
			$\frac{1}{32}$  & $\frac{1}{32}$ & $\frac{1}{32}$ & $\frac{1}{32}$ & $1,0,\frac{1}{4}$  \\
			&&&&\\
			\hline
		\end{tabular} 
		
		\vspace{0.5cm}
		Table 1. Fusion Rules.
	\end{center}	 
	
	\medskip\noindent
	The meaning of the fusion rules is the inclusion 
	$$V_\lm^{(a)} \cdot V_\mu^{(a)} \subseteq \bigoplus_{\nu \in Sp(\lm,\mu)} V_\nu^{(a)}$$
where $\lm,\mu \in Sp$ and $Sp(\mu,\lm)$
is the $(\lm,\mu)$-entry in Table~1.
	
	\medskip
	\noindent
	\begin{defi}
Let $(V,(~,~),\cdot)$ be a triple satisfying $(M1)$ and $(M2)$,
let $A$ be a set of Majorana axes in $V$ satisfying $(M3)$ to $(M7)$,
and let $G$ be the subgroup of the automorphism group of $(V,(~,~),\cdot)$
generated by the Majorana involutions associated with the axes in $A$.
Then the triple $(V,(~,~),\cdot)$ is called a {\em Majorana algebra}, and
the quintuple $(V,(~,~),\cdot,A,G)$ is called 
a {\em Majorana representation} of $G$.
	\end{defi}

\medskip\noindent
If $(V_M,(~,~),\ast)$ is the full name of the Monster algebra
and $2A$ denotes the set of $2A$-axes in this algebra,
then $(V_M,(~,~),\ast,2A,M)$ is a Majorana representation of the Monster.
There is a 1-1 correspondence between $2A$-involutions $t$ in the monster $M$
and $2A$-axes $a$ in $2A$ such that if $\phi \colon M \to GL(V_M)$ is
the given representation of $M$ on $V_M$, then $\phi(t) = \tau(a)$.
It follows that any subgroup $G$ of $M$ generated by a set of $2A$-involutions
possesses at least one Majorana representation.

\medskip\noindent
The above definitions were motivated by a theorem proved by S.~Sakuma in \cite{sak07} making use of earlier results by M.~Miyamoto. In Majorana terms this theorem as stated and proved in \cite{ipss09} sounds as follows.

\begin{thm} \label{sakuma} {\rm (Sakuma)}
	Let $(V,(~,~),\cdot)$ be a Majorana algebra, let $A=\{a_0,a_1\}$ for $a_0 \ne a_1$ be Majorana axes, and let $G$ be the dihedral subgroup $D_{2n}$ generated by $\tau(a_0)$ and $\tau(a_1)$, where $n$ is the order of the product of the generators. Then $n \le 6$ and there are at most eight possibilities for the isomorphism type of the subalgebra generated by $A$ in the Majorana representation $(V,(~,~),\cdot,A,G)$ of $G$. 
\end{thm}

\begin{center} 
	\noindent
	\begin{tabular}{|c|c|c|}
		\hline
		&&\\
		Type & Basis & Products and angles \\
		&&\\
		\hline
		&&\\
		2A & $a_0,a_1,a_\rho$ & $a_0 \cdot a_1=\frac{1}{2^3}(a_0+a_1-a_\rho),~a_0 \cdot a_\rho=\frac{1}{2^3}(a_0+a_\rho-a_1)$ \\
		&&$(a_0,a_1)=(a_0,a_\rho)=(a_1,a_\rho)=\frac{1}{2^3}$\\
		&& \\ 
		2B & $a_0,a_1$ &$a_0 \cdot a_1=0$,~$(a_0,a_1)=0$ \\
		&&\\
		&  &$a_0 \cdot a_1=\frac{1}{2^5}(2a_0+2a_1+a_{-1})-\frac{3^3 \cdot 5}{2^{11}}u_\rho$\\
		3A& $a_{-1},a_0,a_1,$ & $a_0 \cdot u_\rho=\frac{1}{3^2}(2a_0-a_1-a_{-1})+\frac{5}{2^5}u_\rho$~~~~\\
		&$u_\rho$& $u_\rho \cdot u_\rho=u_\rho$\\
		&& $(a_0,a_1)=\frac{13}{2^8}$,~$(a_0,u_\rho)=\frac{1}{2^2}$,~$(u_\rho,u_\rho)=\frac{2^3}{5}$
		\\
		&&\\
		3C & $a_{-1},a_0,a_1$ & $a_0 \cdot a_1=\frac{1}{2^6}(a_0+a_1-a_{-1}),~(a_0,a_1)=\frac{1}{2^6}$ \\
		&&\\ 
		&  & ~$a_0 \cdot a_1=\frac{1}{2^6}(3a_0+3a_1+a_2+a_{-1}-3v_\rho)$\\
		4A & $a_{-1},a_0,a_1,$ & $a_0 \cdot v_\rho=\frac{1}{2^4}(5a_0-2a_1-a_2-2a_{-1}+3v_\rho)$\\
		&$a_2,v_\rho$&~$v_\rho \cdot v_\rho=v_\rho$, ~$a_0 \cdot a_2=0$ \\
		& & $(a_0,a_1)=\frac{1}{2^5},~(a_0,a_2)=0,~(a_0,v_\rho)=\frac{3}{2^3},~(v_\rho,v_\rho)=2$\\
		&&\\
		4B & $a_{-1},a_0,a_1,$ & $a_0 \cdot a_1=\frac{1}{2^6}(a_0+a_1-a_{-1}-a_2+a_{\rho^2})$
		\\
		& $a_2,a_{\rho^2}$ & $a_0 \cdot a_2=\frac{1}{2^3}(a_0+a_2-a_{\rho^2})$\\
		&& $(a_0,a_1)=\frac{1}{2^6},~(a_0,a_2)=(a_0,a_\rho)=\frac{1}{2^3}$ \\
		&&\\
		&& $a_0 \cdot a_1=\frac{1}{2^7}(3a_0+3a_1-a_2-a_{-1}-a_{-2})+w_\rho$
		\\
		5A & $a_{-2},a_{-1},a_0,$ & $a_0 \cdot a_2=\frac{1}{2^7}(3a_0+3a_2-a_1-a_{-1}-a_{-2})-w_\rho$
		\\
		& $a_1,a_2,w_\rho$ & $a_0 \cdot w_\rho=\frac{7}{2^{12}}(a_{1}+a_{-1}-a_2-a_{-2})+\frac{7}{2^5}w_\rho$\\
		& & $w_\rho \cdot w_\rho=\frac{5^2 \cdot 7}{2^{19}}(a_{-2}+a_{-1}+a_0+a_1+a_2)$\\
		&&$(a_0,a_1)=\frac{3}{2^7},~(a_0,w_\rho)=0$, $(w_\rho,w_\rho)=\frac{5^3 \cdot 7}{2^{19}}$\\
		&& \\
		& & $a_0 \cdot a_1=\frac{1}{2^6}(a_0+a_1-a_{-2}-a_{-1}-a_2-a_3+a_{\rho^3})+\frac{3^2 \cdot 5}{2^{11}}u_{\rho^2}$\\
		6A& $a_{-2},a_{-1},a_0,$ &$a_0 \cdot a_2=\frac{1}{2^5}(2a_0+2a_2+a_{-2})-\frac{3^3 \cdot 5}{2^{11}}u_{\rho^2}$  \\ 
		&$a_1,a_2,a_3$  &$a_0 \cdot u_{\rho^2}=\frac{1}{3^2}(2a_0-a_2-a_{-2})+\frac{5}{2^5}u_{\rho^2}$  \\
		&$a_{\rho^3},u_{\rho^2}$ & $a_0 \cdot a_3=\frac{1}{2^3}(a_0+a_3-a_{\rho^3})$, $a_{\rho^3} \cdot u_{\rho^2}=0$, $(a_{\rho^3},u_{\rho^2})=0$\\
		&&$(a_0,a_1)=\frac{5}{2^8}$, $(a_0,a_2)=\frac{13}{2^8}$, $(a_0,a_3)=\frac{1}{2^3}$\\
		&&\\
		\hline
	\end{tabular} 
	
	\vspace{0.5cm} 
	Table 2. Norton--Sakuma Algebras.
\end{center}

\medskip\noindent
Each of the eight algebras in Sakuma's theorem are subalgebras in the Monster algebra, their explicit forms were computed by S.~P.~Norton in \cite{nor96} and given in Table~2 with respect to the Majorana scaling.
The name of the algebra generated by $a_0$ and $a_1$ is the conjugacy class of the Monster containing the product $\tau(a_0)\tau(a_1)$.

\medskip\noindent
Therefore, Sakuma's theorem provides the classification of Majorana representations of the dihedral groups. On progress in classifying Majorana representations of further groups we refer the reader to surveys \cite{iv18}, \cite{iv22}.

\medskip\noindent
Our main result is the following.

\begin{main}
The group $U_3(5)$ has a unique Majorana representation
satisfying {\rm (M8)} below.
This representation has dimension $798$, it is $2$-closed,
is spanned by Majorana and $3A$-axes and is based on
an embedding into the Monster.
\end{main}

\medskip\noindent
Here the {\em dimension} of a Majorana representation $(V,(~,~),\cdot,A,G)$
is the vector space dimension of the subalgebra $\agen{A}$ of $V$
generated by $A$.
The representation is called $i$-{\em closed}, when $\agen{A}$
is spanned by products of at most $i$ elements of $A$.

\medskip\noindent
A $3A$-{\em axis} is a vector with the role of $u_\rho$ in a type $3A$
subalgebra (an idempotent of squared length $\frac85$, expressed in
the Majorana axes $a_{-1},a_0,a_1$ by the first equation in the $3A$ part
of Table~2).
Similarly, $4A$-axes and $5A$-axes are vectors with the role of
$v_\rho$ and $w_\rho$, respectively, in type $4A$ ($5A$) subalgebras.

\medskip\noindent
Condition (M8) is
\begin{itemize}
	\item[{\rm (M8)}]
	The vectors $a_\rho$, $a_{\rho^2}$ and $a_{\rho^3}$
in type $2A$, $4B$ and $6A$ algebras, respectively, are Majorana axes.
The vectors $u_\rho$, $v_\rho$, $w_\rho$ in algebras of type $3A$, $4A$ and
$5A$ depend solely on the group element $\rho = \tau(a_0) \tau(a_1)$
(rather than on the whole dihedral group $\ggen{\tau(a_0),\tau(a_1)}$).
\end{itemize}

\medskip\noindent
Condition (M8) is satisfied by $(V_M,(~,~),\ast,2A,M)$ by \cite{nor96}.
In our $U_3(5)$ setting this condition will allow us to identify $3A$-axes
with subgroups of order 3.
Note that condition (M8) implies conditions (2A) and (3A) in \cite{ipss09}.

\section{On the $U_3(5)$ group}	
In this section we summarise the required properties of the group $U_3(5)$. They can be deduced from the information in the relevant section of the ATLAS \cite{ccnpw} and the properties of the Hoffman--Singleton graph found in \cite{bm22}, \S10.19.

\begin{itemize}
	\item[{\rm (P1)}] $U_3(5)$ is a simple group of order $126\,000=2^4\cdot 3^2\cdot 5^3 \cdot 7$ with a unique class of involutions. Its Schur multiplier is of order 3 and the outer automorphism group is isomorphic to $S_3$. 
	\item[{\rm (P2)}] $U_3(5)$ contains three classes of subgroups isomorphic to $A_7$, which are maximal and the classes are transitively permuted by the outer automorphisms.
	\item[{\rm (P3)}] The action of $U_3(5)$ on the cosets of a subgroup $A_7$ preserves a structure of the Hoffman--Singleton graph, which is the unique strongly regular graph with parameters $v=50$, $k=7$, $\lambda=0$, $\mu=1$. The action of $U_3(5)$ on this graph has rank 3.
	\item[{\rm (P4)}] The action of an $A_7$-subgroup on the cosets of an $A_7$ subgroup from a different class has two orbits with lengths 15 and 35 and stabilizers $L_3(2)$ and $(S_3 \times S_4)^{+}$, the latter being the normalizer of a subgroup of order 3 in both $A_7$'s.
	\item[{\rm (P5)}] If $t$ is an involution in $U_3(5)$, then the fixed vertices of $t$ on the Hoffman--Singleton graph form a Petersen subgraph on which the centralizer $C(t)$ of $t$ in $U_3(5)$ induces an action of $S_5$ with kernel $\ggen{t}$. The whole $C(t)$ is a non-split extension of $\ggen{t}$ by $S_5$ in which a transposition of $S_5$ lifts to an involution. $C(t)$ is the unique index 2 subgroup of $GL_2(5)$ with the class of 20 non-central involutions. $C(t)$ is maximal in $U_3(5)$. 
	
	\item[{\rm (P6)}] The action of $U_3(5)$ on its 525 involutions
has rank 8, with suborbits of sizes 1, 20, 120, 120, 120, 48, 48, 48,
corresponding to product orders 1, 2, 3, 4, 6, 5, 5, 5, respectively.
The action of $U_3(5).S_3$ has rank 6, the three suborbits of size 48 fuse
into a single suborbit of size 144.

\item[{\rm (P7)}] The Monster group contains a maximal 5-local subgroup
$$N(5A^2) \cong (5^2:4 \cdot 2^2 \times U_3(5)):S_3$$ 
where
$$N(5A) \cong (D_{10} \times HN)\cdot 2$$
(see \cite{w88}).

\item[{\rm (P8)}] The involutions in the $U_3(5)$-subgroup in $N(5A^2)$ are $2A$-involutions in the Monster by Table~5 in \cite{n98}. 

\end{itemize}

\section{An upper bound}
Property (P8) implies that $U_3(5)$ possesses a Majorana representation based
on embedding into the Monster. It would be very hard 
to deduce any exact information on this representation, like the dimension,
by pre-Majorana methods. The characterization \cite{fim16} of the Majorana representation of the Harada--Norton group is not explicit enough either. What one can do is to get an upper bound on the dimension of this particular representation. The procedure is rather standard.

\begin{lem} \label{exact}
The following assertions hold:
\begin{itemize}
	\item[{\rm (i)}] $N_M(U_3(5))/U_3(5) \cong 5^2:(SL_2(3) * Z_4)$;
	\item[{\rm (ii)}] $C_M(U_3(5)) \cong 5^2:(Q_8*Z_4)$. 
\end{itemize}
\end{lem}	

\pf
Part (i) follows from (P7), and (ii) is immediate from (i).~\qed

\medskip\noindent
It is clear that the Majorana subalgebra of $U_3(5)$ in the Monster algebra
is contained in 
$$C_{V_M}(C_M(U_3(5))).$$
If $b$ is the dimension of the above centralizer, then the dimension
of the Majorana representation of $U_3(5)$ in $V_M$ is at most $b$.
On the other hand, $b$ can be calculated by restricting the character
of $M$ on $V_M$ to $C_M(U_3(5))$
and counting the number of trivial components. 

\medskip\noindent
The character table of $C_M(U_3(5))$ was found on the internet
by William Giuliano to whom we are very thankful.
The fusion of classes was recovered via the embedding 
$$C_M(U_3(5)) \le C_M(A_5) \cong A_{12}.$$

\begin{center}
\begin{tabular}{|c|c||c|c||c|c||c|c|}
	\hline
	$A_{12}$&$M$&$A_{12}$&$M$&$A_{12}$&$M$&$A_{12}$&$M$\\
	\hline
	$(1^{12})$&$1A$&$(4,2)$&$4B$&$(6,2)$&$6C$&$(4,3,2)$&$12C$\\
	\hline
	$(2^2)$&$2A $&$(4,2^3) $&$4B $&$(7) $&$7A $&$(4,3^2,2) $&$12C $\\
	\hline
	$(2^6) $&$2A $&$(5) $&$5A $&$(8,2) $&$8B $&$(6,4) $&$12C $\\
	\hline
	$(2^4) $&$2B $&$(5^2) $&$5A $&$(8,4) $&$8B $&$(7,2^2) $&$14A $\\
	\hline
	$(3) $&$3A $&$(3,2^2) $&$6A $&$(9) $&$9A $&$(5,3) $&$15A $\\
	\hline
	$(3^2) $&$3A $&$(3^2,2^2) $&$6A $&$(9,3) $&$9A $&$(5,3^2) $&$15A $\\
	\hline
	$(3^4) $&$3A $&$(6,2^3) $&$6A $&$(5,2^2) $&$10A $&$(5,4,2) $&$20B $\\
	\hline
	$(3^3) $&$3B $&$(6^2) $&$6A $&$(10,2) $&$10A $&$(7,3) $&$21A $\\
	\hline
	$(4^2) $&$4A $&$(6,3,2) $&$6B $&$(11) $&$11A $&$(5,3,2^2) $&$30B $\\
	\hline
	$(4^2,2^2) $&$4A $&$(3,2^4) $&$6C $&$ (4^2,3)$&$12A $&$(7,5) $&$35A $\\
	\hline
\end{tabular}

\vspace{0.5cm}
Table 3. Fusion of $A_{12}$-classes into the Monster.
\end{center}

\medskip\noindent
The fusion of $A_{12}$-classes into the Monster is well known,
appeared for instance in \cite{l17} and we present it here in Table~3
for future reference.  
 
\medskip\noindent
Since $5^2$ is a Sylow 5-subgroup of $A_{12}$ and $N_{A_{12}}(5^2)$ is the intersection with $A_{12}$ of $N_{S_{12}}(5^2)\cong (F_{20} \wr 2) \times S_2$ and $C_M(U_3(5))$ as in Lemma~\ref{exact}(ii) is rather visible inside this normalizer, we obtain the final result of these calculations.
\begin{lem}
	$$\dim C_{V_M}(C_M(U_3(5)))=990.$$
\qed
\end{lem} 

\section{Subrepresentations and relations}
We start this section by presenting the list of eigenvectors of a Majorana axis $a_0$ inside Norton--Sakuma algebras. This list is well known \cite{ipss09} and can be deduced from the product rules.
In order to save space, we write $e_1 := a_1 - a_{-1}$
and $e_2 := a_2 - a_{-2}$.

\begin{center}
	\begin{tabular}{@{~}c@{~~}|c|c@{~~}|@{~~}c@{~}}
		\hline
		Type & 0 & $\frac{1}{4}$ & $\frac{1}{32}$\rule{0pt}{11pt}\\[2pt]
		\hline
		$2A$    & $a_1+a_\rho-\frac{1}{2^2}a_0$   &$a_1-a_\rho$    &   \\
		$2B$    &  $a_1$ & & \\
		$3A$    &  $u_\rho-\frac{2 \cdot 5}{3^3}a_0+\frac{2^5}{3^3}(a_1+a_{-1})$& $u_\rho-\frac{2^3}{3^2 \cdot 5}a_0-\frac{2^5}{3^2 \cdot 5}(a_1+a_{-1})$
		& $e_1$\\
		$3C$ & $a_1+a_{-1}-\frac{1}{2^5}a_0$ &  & $e_1$ \\
		$4A$ & $v_\rho-\frac{1}{2}a_0+2(a_1+a_{-1})+a_2$,~~~$a_2$ & $v_\rho-\frac{1}{3}a_0-\frac{2}{3}(a_1+a_{-1})-\frac{1}{3}a_2$ & $e_1$ \\
		$4B$ & $a_1+a_{-1}-\frac{1}{2^5}a_0-\frac{1}{2^3}(a_{\rho^2}-a_2)$, & $a_2-a_{\rho^2}$ & $e_1$ \\[1pt]
		& $a_2+a_{\rho^2}-\frac{1}{2^2}a_0$ & & \\[1pt]
		$5A$ & $w_\rho+\frac{3}{2^9}a_0-\frac{3 \cdot 5}{2^7}(a_1+a_{-1})-\frac{1}{2^7}(a_2+a_{-2}),$ & $w_\rho+\frac{1}{2^7}(a_1+a_{-1}-a_2-a_{-2})$ & $e_1,$\\[1pt]
		&$w_\rho-\frac{3}{2^9}a_0+\frac{1}{2^7}(a_1+a_{-1})+\frac{3 \cdot 5}{2^7}(a_2+a_{-2})$ & & $e_2$ \\
		$6A$ & $u_{\rho^2}+\frac{2}{3^2 \cdot 5}a_0-\frac{2^4}{3^2 \cdot 5}
		(a_1+a_{-1})-~~~$ & $u_{\rho^2}-\frac{2^3}{3^2 \cdot 5}a_0-$\hspace{75pt} & $e_1,$ \\[1pt]
		&$~~~~~\frac{2^3}{3^2 \cdot 5}(a_2+a_{-2}+a_3-a_{\rho^3}),$&$~~~~~\frac{2^3}{3^2 \cdot 5}(a_2+a_{-2}+a_3-a_{\rho^3}),$& \\[1pt]
		&$a_3+a_{\rho^3}-\frac{1}{2^2}a_0$, $u_{\rho^2}-\frac{2 \cdot 5}{3^3}a_0+\frac{2^5}{3^3}(a_2+a_{-2})$ & $a_3-a_{\rho^3}$ & $e_2$ \\
		&&&\\
		\hline
		
	\end{tabular}	

	\vspace{0.5cm}
	Table 4. Eigenvectors of Norton--Sakuma algebras.
\end{center}

\medskip\noindent
Next we formulate and prove an important special case of the {\em resurrection principle} \cite{ipss09}.
In what follows $\agen{X}$ denotes the subalgebra generated by a set $X$.
We shall write $ab$ instead of $a \cdot b$.

\begin{lem}
	Let $(V,(~,~),\cdot,A,G)$ be a Majorana representation of a group $G$. Suppose that $b_0,b_1,b_2$ are Majorana axes such that 
	$$B_1:=\agen{b_0,b_1} {\rm ~and~}B_2:=\agen{b_0,b_2}$$
	are $3A$-algebras with $3A$-axes $u_1$ and $u_2$, respectively. Then the product $u_1 u_2$ is explicitly expressible in terms of products of Majorana axes and products of Majorana axes with a $3A$-axis, and similarly for the inner product $(u_1,u_2)$.
\end{lem}

\pf Let $\alpha_1$ be the $0$-eigenvector of ${\rm ad}(b_0)$ in $B_1$,
$\alpha_2$ and $\beta_2$ the $0$- and $\frac{1}{4}$-eigenvectors of
that operator in $B_2$ normalized as in Table~4. Then by the fusion rule
$\alpha_1 \alpha_2$ is a 0-eigenvector, while $\alpha_1 \beta_2$ is a
$\frac{1}{4}$-eigenvector of ${\rm ad}(b_0)$. These eigenvectors
have $u_1 u_2$ as the leading terms and the remaining terms
are products of two Majorana axes and of a Majorana axis with $u_1$ or $u_2$.
Now $u_1 u_2$ can be found from
$$b_0 (\alpha_1 \alpha_2-\alpha_1 \beta_2)=-\frac{1}{4}\alpha_1 \beta_2=-\frac{1}{4}u_1 u_2+\cdots .$$
For the inner product we just apply $0=(\alpha_1,\beta_2)=(u_1,u_2)+\cdots$.
\qed

\begin{lem} \label{rel}
	Let $G$ be the semidirect product of an elementary abelian group
of order $9$ generated by elements $x$ and $y$, and a group of order $2$
generated by $t$ which inverts both $x$ and $y$. Then $G$ possesses
a unique Majorana representation satisfying {\rm (M8)}
such that any two Majorana axes generate a $3A$ subalgebra.
Furthermore, if $b_1,\ldots,b_9$ and $u_1,\ldots,u_4$ are the Majorana
and $3A$-axes in this representation, then
	$$45 \sum_{i=1}^4u_i-32\sum_{j=1}^9 b_j=0$$
\end{lem}
\medskip
\noindent
The relations in Lemma~\ref{rel} were called in \cite{iv11}
{\em Pasechnik relations}  and this terminology became standard,
although the referee of \cite{iv11} pointed out that these relations
were known for a long time in the theory of Vertex Operator Algebras.

\medskip\noindent
The Majorana representation of $A_7$ based on the embedding into the Monster
was characterised in \cite{iv11}.

\begin{lem} \label{a7}
	The group $A_7$ possesses a unique Majorana representation
satisfying {\rm (M8)}. This representation has dimension $196$.
It is spanned by $105$ Majorana axes, which are linearly independent
and $140$ $3A$-axes associated with order $3$ elements of cyclic type $3^2$.
The latter are subject to $49$ linearly independent relations
of which $35$ are Pasechnik relations and $14$ are Faradzev relations.
All the $5A$-axes in the representation are linear combinations of
the Majorana and $3A$-axes.\qed  
\end{lem}

\medskip\noindent
An explicit form of Faradzev relation was obtained by Clara Franchi and Mario Mainardis \cite{fm21} as an alternating sum of some 48 $3A$-axes. 
This form played a crucial role in some early stages of our project.

\medskip\noindent
The information in the following lemma was obtained by J.~McInroy by the expansion algorithm for Majorana algebra described in \cite{ms20}. We are extremely thankful to Justin for his sharing this information with us.

\begin{lem} \label{2S5}
	The group $C=2\cdot S_5$, which is the centralizer of an involution in $U_3(5)$, possesses a unique Majorana representation. This representation has dimension $31$, and is spanned by twenty-one Majorana and ten $3A$-axes. Furthermore, if $\rho$ is an element of order $3$ in $C$, then 
	\begin{itemize}
		\item[{\rm (i)}] for an involution $t$ in $C$ which generates with $\rho$ a $GL_2(3)$-subgroup, we have $(u_\rho,a_t)=\frac{1}{36}$;
		\item[{\rm (ii)}] for an element $\sigma$ of order order $3$ in $C$ which generates with $\rho$ an $SL_2(3)$-subgroup, we have  $(u_\rho,u_\sigma)=\frac{8}{81}$;
		\item[{\rm (iii)}] for an element $\pi$ of order $3$ in $C$, which generates with $\rho$ an $SL_2(5)$-subgroup, we have $(u_\rho,u_\pi)=\frac{16}{405}$;
	\end{itemize}	
	\qed	 
\end{lem}

\section{Shape and $1$-closure}
We start our construction of a Majorana representation of $G=U_3(5)$ by setting a vector space with a set $A$ of $525$ vectors $a_t$ indexed by the involutions $t \in T$, where $T$ is the class of involutions in $G$. In order to proceed we need to find the {\em shape} of the representation we are aiming at, which is a map $sh$ from the set of dihedral subgroups in $G$ into the set of Norton--Sakuma algebras such that 
$$\agen{a_t,a_s} \cong sh(\ggen{t,s}_G).$$
\begin{lem} \label{shape}
	The shape $sh$ of a representation of $G \cong U_3(5)$ is uniquely determined and the set of images of $sh$ is 
	$$\{2A,3A,4B,5A,6A\}$$
\end{lem}

\pf There is only one algebra for $D_{12}$ and one for $D_{10}$. By (P6) every subgroup $D_6$ is contained in a $D_{12}$-subgroup. Since a $6A$ algebra contains $3A$- but not $3C$-subalgebras, the image of a $D_6$-subgroup is $3A$. Similarly, since $6A$ contains $2A$- but not $2B$-subalgebras, the image of a $D_4$-subgroup is $2A$. Finally $4A$ contains a $2B$-subalgebra, which we do not have, therefore, the image of $D_8$ is $4B$.~\qed

\medskip
\noindent
Lemma~\ref{shape} enables us to determine the form on the linear span of $A$. 

\begin{cor} \label{innerprod22}
	The inner product $(a_t,a_s)$ is equal to $1$, $\frac{1}{8}$, $\frac{13}{256}$, $\frac{1}{64}$, $\frac{3}{128}$ and $\frac{5}{256}$ when $a(st)$ is equal to $1$, $2$, $3$, $4$, $5$ and $6$, respectively.~\qed
\end{cor}

\medskip\noindent
Computationally we obtain the following.

\begin{lem}
	The $525 \times 525$ Gram matrix $\Gamma(A)=||(a_t,a_s)||$ has full rank.~\qed
\end{lem}

\medskip
\noindent
Thus the 1-closure, which is the subspace spanned by the Majorana generators of a representation of $U_3(5)$, is $525$-dimensional with uniquely determined inner product on it.

\section{$2$-closure}
Next we consider the 2-closure, which is the space generated by the Majorana generators together with their pairwise products. By Lemma~\ref{shape} the 2-closure is spanned by the Majorana generators together with $3A$- and $5A$-axes. The following lemma takes care of the $5A$-axes.

\begin{lem} \label{5A}
	Every $5A$-axis of a Majorana representation of $U_3(5)$ is a linear combination of the Majorana generators and $3A$-axes in the $2$-closure.
\end{lem}

\pf The result follows from Lemma~\ref{a7} since there are three $U_3(5)$-orbits on $D_{10}$-subgroups by (P6) and three orbits on the $A_7$-subgroups by (P2) with respective inclusion.~\qed

\begin{lem}
	The $2$-closure of a Majorana representation of $U_3(5)$ is spanned by the $525$ Majorana generators and $1750$ $3A$-axes indexed by the subgroups of order $3$ in $U_3(5)$.
\end{lem}

\pf By Lemmas~\ref{shape} and \ref{5A} we should only consider $3A$-axes. By (P6) there are
$$525\cdot 120/6=10500$$
dihedral groups of order $6$ in $U_3(5)$, six for every subgroup of order 3.
Comparing the orders of normalizers we conclude that the normalizer of
a $3$-subgroup is contained in some $A_7$-subgroup. So the glueing of
$3A$-subalgebras already takes place in $A_7$-subrepresentations
(cf. Lemma~\ref{a7}).~\qed

\medskip\noindent
We denote by 
$$U=\{u_\rho \mid \rho \in U_3(5), \rho^3=1\}$$
the set of $3A$-axes in the representation, understanding that $u_\rho=u_{\rho^{-1}}$.

\medskip
\noindent
Next we are aiming to determine the Gram matrix of the set $A \cup U$ and by calculating its rank decide on the dimension of the $2$-closure.
The following result was obtained computationally.

\begin{lem} \label{innerprod23}
The pairs of Majorana $2A$- and $3A$-axes in a representation of $U_3(5)$
are as described in Table~$5$, where $T_H$ is the number of involutions generating a group $H$ of a given isomorphism type together with a fixed order $3$-subgroup generated by $\rho$.~\qed
\end{lem}

\begin{center}
 	\begin{tabular}{r|c|c|c}
 		$|T_H|$ & $(t\rho)^M$ &$H = \ggen{t,\rho}$ & $(a_t,u_\rho)$ \\[1pt]
 		\hline
 		$3$ & $6A$& $C_6$ & 0 \\[1pt]
 		$18$ & $2A$ & $S_3$ &$\frac{1}{4}$\\[1pt]
 		$36$ & $3A$ & $A_4$ &$\frac{1}{9}$ \\[1pt]
 		$108$ & $4B$ & $S_4$ &$\frac{1}{36}$\\[1pt]
 		$36$ & $8C$ & $GL_2(3)$ &$\frac{1}{36}$\\[1pt]
 		$108$ & $5A$ & $A_5$ & $\frac{1}{18}$\\[1pt]
 		$216$ & $7A$ & $L_2(7)$ & $\frac{1}{24}$\\[1pt]
 		\hline
 		$525$&&&\\
 	\end{tabular}

\vspace{0.5cm}
Table 5. $(2A,3A)$-pairs in $U_3(5)$
\end{center}

\medskip\noindent
Notice that all $(2A,3A)$-pairs inside the Monster were classified by Simon Norton and are presented in Table~3 of \cite{nor96} (in a scaling different from ours). We do not have a Majorana version of this classification and can only compare the results and be happy when they are consistent, which is the case here. The entries in the second column showing the class of $(t\rho)$ in the Monster are taken from Norton's table. The inner products can be calculated inside $A_7$-subrepresentations as in Lemma~\ref{a7} for all pairs except those generating $GL_2(3)$, where we apply Lemma~\ref{2S5}. 

\medskip
\noindent
The inner products between $3A$-axes were also determined computationally.

\begin{lem} \label{innerprod33}
	Let $\rho$ be a subgroup of order $3$ in $U_3(5)$. Then all the order $3$ subgroups $\sigma$ in $U_3(5)$ are classified by the isomorphism type of the subgroup $H=\ggen{\rho,\sigma}$ generated by $\rho$ and $\sigma$ as indicated in Table~$6$.
\end{lem}

\begin{center}
	\begin{tabular}{r|r|c|c}
		$|U_H|$ & $|H|$ &$H=\ggen{\rho,\sigma}$ & $(u_\rho,u_\sigma)$ \\[2pt]
		\hline
		$1$ & $3$ & $C_3$ & $\frac{8}{5}$ \\[2pt]
		$12$ & $9$ & $C_3 \times C_3$ & $0$ \\[2pt]
		$36$ & $12$ & $A_4$ & $\frac{136}{405}$\\[2pt]
		$144$ & $21$ & $F_7^3$ & $\frac{4}{27}$\\[2pt]
		$18$ & $24$ & $SL_2(3)$ & $\frac{8}{81}$\\[2pt]
		$72$ & $36$ & $C_3 \times A_4$ & $\frac{64}{405}$ \\[2pt]
		$54$ & $60$ & $A_5$ & $\frac{16}{405}$\\[2pt]
		$9$ & $120$ & $SL_2(5)$ & $\frac{16}{405}$ \\[2pt]
		$108$ & $168$ & $L_2(7)$ & $\frac{32}{405}$ \\[2pt]
		$216$ & $168$ & $L_2(7)$ & $\frac{4}{81}$\\[2pt]
		$216$ & $360$ & $A_6$ & $\frac{32}{405}$\\[2pt]
		$216$ & $2520$ & $A_7$ & $\frac{8}{81}$ \\[2pt]
		$216$ & $2520$ & $A_7$ & $\frac{32}{405}$ \\[2pt]
		$432$ & $126000$ & $U_3(5)$ & $x$\\[2pt]
		\hline
		1750 &&&\\		
	\end{tabular}

\vspace{0.5cm}
Table 6. Inner products of $3A$-axes in $U_3(5)$.
\end{center}

\medskip\noindent
In Table~6 $U_H$ is the set of subgroups $\sigma$ generating with $\rho$
a subgroup isomorphic to $H$ with further subdivisions in the two cases
$L_2(7)$ and $A_7$. For $L_2(7)$ the inner product is $\frac{32}{405}$
when $\rho$ and $\sigma$ are normalized by a common involution in $H$,
and $\frac{4}{81}$ otherwise. For $A_7$ the inner product is $\frac{8}{81}$
when the generating order 3 subgroups in $A_7$ are both of cycle type $3^2$,
and $\frac{32}{405}$ when the generation is by a 3- and a $3^2$-element.

\medskip\noindent
The inner products are calculated in $A_7$ via Lemma~\ref{a7}
with \cite{iv11} and in $C(t) \cong 2 \cdot S_5$ as in  Lemma~\ref{2S5}.
In this way, proceeding inductively, we cannot calculate the inner product
in the case when $\rho$ and $\sigma$ generate the whole group $G$
so we put $x$ in the relevant position.
The reason of having just one $x$ is justified by the following result
obtained computationally.

\begin{lem}
	Let $\rho$ be a subgroup of order $3$ in $G=U_3(5)$ and let $\Gamma_{-}(\rho)$ denote the set of order $3$ subgroups in $U_3(5)$ which generate with $\rho$ the whole group $G$. Then
	\begin{itemize}
		\item[{\rm (i)}] $\Gamma_{-}(\rho)$ is a union of six regular orbits of $N:=N_{G}(\ggen{\rho}) \cong (S_3 \times S_4)^+$ of order $72$;
		\item[{\rm (ii)}] the orbits in $(i)$ are transitively permuted by 
		$N_{{\rm Aut}(G)}(\ggen{\rho})/N \cong S_3$;
		\item[{\rm (iii)}] if $\sigma \in \Gamma_{-}(\rho)$, then there is no involution in $G$ which normalizes both $\ggen{\rho}$ and $\ggen{\sigma}$.
	\end{itemize}	
\end{lem}

\medskip\noindent
The following lemma was obtained computationally. We sketch the setting.

\medskip\noindent
Let $V^{(2)}:=\lspan{A\cup U}$, the linear span of $A \cup U$, be
the 2-closure of the representation, and for a $3A$-axis $u \in U$
corresponding to a subgroup $\rho$ of order 3, let $V_+(u)$ be
the subspace of $V^{(2)}$ spanned by all Majorana axes in $A$
together with all the $3A$-axes except those in the last row of Table~5
(that is, except for those corresponding to a subgroup $\sigma$ of order 3
with $\ggen{\rho,\sigma} = U_3(5)$).
Notice that for every $v \in V_+(u)$ the product $u v$ is contained
in an $A_7$- or $2.S_5$ subrepresentation, thus can be computed and is
contained in $V^{(2)}$. Let $\Gamma_{-}(u)$ be the set of $3A$-axes
corresponding to the last row in Table~5, so that $|\Gamma_{-}(u)|=432$. Let 
$$\Gamma_{-}(u) = O_1 \cup O_2 \cup \dots \cup O_6$$
be the decomposition of $\Gamma_{-}(u)$ into a disjoint union
of $N_G(\ggen{\rho})$-orbits. For $1 \le i \le 6$, let $F_i$
be the function on $V^{(2)}$ which is 1 on $O_i$ and 0 on
the remaining $O_j$'s.  

\begin{lem}
	The following assertions hold:
\begin{itemize}
	\item[{\rm (i)}] $x=\frac{4}{81}$;
	\item[{\rm (ii)}] $\dim V^{(2)}=798$;
	\item[{\rm (iii)}] $\dim V_+(u)=796$;
	\item[{\rm (iv)}] there is an orthogonal complement $V_-(u)$ of $V_+(u)$ in $V^{(2)}$ spanned by
	$$F_1+F_2-F_3-F_4,~~~-F_1-F_2+F_5+F_6,~~~F_3+F_4-F_5-F_6$$
	for some arrangement of the orbits $O_i$ into pairs.
\end{itemize}
\end{lem}

\pf The left hand side of the relations from an $A_7$-subrepresentation are vectors of length zero and they must be zero vectors since the form is positive-definite. In particular they have to be perpendicular to all other vectors, which enabled us to determine $x$, proving (i). The power space  of $A \cup U$ factorized over the linear span of the relations in the three classes $A_7$-subrepresentations turned out to be $796$-dimensional positive-definite which gives (ii). The assertions (iii) and (iv) were also achieved computationally.~\qed

\medskip\noindent
It is clear from the above lemma that in order to close the product on $V^{(2)}$ it suffices to show that for a triple of $3A$-axes $u_1$, $u_2$, $u_3$ such that $u_2$ and $u_3$ are in the same regular orbit of $G(u_1)$ on $\Gamma_{-}(u_1)$, the product
$$u_1 (u_2+u_3)$$   
belongs to $V^{(2)}$. This was achieved by a version of the resurrection principle.

\medskip\noindent
Computationally the following was established.

\begin{lem}
	Let $u_1 \in U$ correspond to an element $\rho$ and $u_2 \in \Gamma_-(u_1)$ to an element $\sigma$. Then there exists an involution $t$ in $G$, such that $\ggen{\rho,t} \cong S_3$ and $\ggen{\sigma,t} \cong S_4$.
\end{lem}

\medskip\noindent
By the above lemma and the shape of the representation of $G$ we have that
the subalgebra $A_1$ generated by $u_1$ and $a(t)$ is the $3A$-algebra and
the subalgebra $A_2$ generated by $u_2$ and $a(t)$ is the 13-dimensional
$S_4$-algebra of shape $(2A,3A)$ as described in \cite{ipss09}. 

\medskip\noindent
% We adopt notation in that paper p.2460 so that ...
We assume that in $A_1$ $a(t)=a_0$ and $u_1=u_\rho$ while in $A_2$ $a(t)=a_{(ij)}$ and $u_2=u_i$. Then $u_j$, which is the image of $u_i$ under $t$ is contained with $u_i$ in the same $G(u_1)$-orbit on $\Gamma_-(u_1)$, since $t \in G(u_1)$. The subgroups of the axes $u_k$ and $u_l$ and normalized by $t$ and therefore, they are in $V_+(u_1)$.
Then the product rule in $A_2$ on p.\,2460 in \cite{ipss09} shows that
$$\alpha_1:=u_i+u_j-\frac{1}{8}(u_k-\frac{8}{45}a_{(ij)}-
\frac{32}{45}(a_{(il)}+a_{(jl)}))$$
$$-\frac{1}{8}(u_l-\frac{8}{45}a_{(ij)}-
\frac{32}{45}(a_{(ik)}+a_{(jk)}))$$
$$-\frac{1}{18}a_{(ij)}-\frac{8}{45}(a_{(kl)}-a_{(ij)(kl)})$$
is a 0-eigenvector of $a_{(ij)}$.
Notice that the expressions in first, second and third brackets are
$\frac{1}{4}$-eigenvectors of $a_{(ij)}$ in the algebras $\agen{a_{(ij)},u_k}$,
$\agen{a_{(ij)},u_l}$ and $\agen{a_{(ij)},a_{(kl)}}$, respectively,
where these algebras are of types $3A$, $3A$ and $2A$, respectively.

\medskip\noindent
Let $\alpha_2$ and $\beta_2$ be $0$- and $\frac{1}{4}$-eigenvectors
of $a(t) = a_0$ in the algebra $A_1 \cong 3A$ as in Table~4. Then by
the fusion rule
$$\alpha_1\alpha_2=u_1(u_i+u_j)+v_1$$
is a 0-eigenvector of $a(t)$ for some $v_1 \in V_{+}(u_1)$, while 
$$\alpha_1\beta_2=u_1(u_i+u_j)+v_2$$
is a $\frac{1}{4}$-eigenvector of $a(t)$, where $v_2 \in V_+(u_1)$.
%
% In fact, $u_k$, $u_l$ are in $V_+(u_1)$ since the three $3A$-axes
% are inverted by $t$.
%

\medskip
\noindent
Now we apply the {\em resurrection principle}, Lemma 1.8 in \cite{ipss09}, compare Lemma 5.

\medskip\noindent
Consider
$$\alpha_1\alpha_2-\alpha_1\beta_2=v_3$$
where $v_3$ is in $V_+(u_1)$ so we can explicitly calculate 
$$a_{(ij)}(\alpha_1\alpha_2-\alpha_1\beta_2)=a_{(ij)} v_3,$$
although by the eigenvalue properties the above expression is equal to
$$0(\alpha_1\alpha_2)-\frac{1}{4}(\alpha_1\beta_2)=-\frac{1}{4}(u_1(u_i+u_j))-\frac{1}{4}v_2,$$
which gives the required expression.

\medskip\noindent
Now we know that the product is closed on $V^{(2)}$, so the latter is the whole Majorana algebra supporting the representation, which completes the proof of Theorem 1.

\medskip\noindent
Since $\alpha_1$ and $\beta_2$ are eigenvectors of $a(t)$ with different
eigenvalues, they are perpendicular. If we expand the equality
$(\alpha_1,\beta_2)=0$ in terms of the above expressions and substitute
the numerical values from Corollary \ref{innerprod22},
Lemma \ref{innerprod23} and Lemma \ref{innerprod33}, we deduce that 
$$(u_1,u_2)=(u_1,u_3)=\frac{4}{81},$$
thus obtaining an independent confirmation of Lemma 17 (i).

\begin{center}
	{\bf Acknowledgement}
\end{center}
We are extremely thankful to our colleagues without whom the project could not be completed: William Giuliano for finding the character table of $C_M(U_3(5))$; to Clara Franchi and Mario Mainardis for deducing the explicit form of the Faradzev relations in the $A_7$-algebra; to the three of them for most careful proofreading of the draft; and to Justin McInroy for sharing with us information on the representation of $2\cdot S_5$.

\end{document}